\documentclass[11pt]{amsart}
\usepackage{amssymb}
\usepackage{amsmath}
\usepackage{amscd}
\usepackage{amsfonts}
\usepackage{mathrsfs}
\usepackage{stmaryrd}
\usepackage{tikz}
\usetikzlibrary{arrows}
\usepackage[T1]{fontenc}
\usepackage{mathtools}
\usepackage{hyperref}

\usepackage{bbm}

\newtheorem{prop-def}{Proposition-Definition}[section]

\begin{document}

\setlength{\oddsidemargin}{0cm} \setlength{\evensidemargin}{0cm}
\baselineskip=18pt

\begin{center}
{\bf{Enumeration of words that contain the pattern $123$ exactly once}}~\\
Mingjia Yang
\end{center}

\subsection*{Abstract:} Enumeration problems related to words avoiding patterns as well as permutations that contain the pattern $123$ exactly once have been studied in great detail. However, the problem of enumerating words that contain the pattern $123$ exactly once is new and will be the focus of this paper. Previously, Doron Zeilberger provided a shortened version of Alexander Burstein's combinatorial proof of John Noonan's theorem that the number of permutations with exactly one $321$ pattern is equal to $\frac{3}{n} \binom{2n}{n+3}$. Surprisingly, a similar method can be directly adapted to words. We are able to use this method to find a formula enumerating the words with exactly one $123$ pattern. Further inspired by Nathaniel Shar and Zeilberger's paper on generating functions enumerating 123-avoiding words with $r$ occurrences of each letter, we examine the algebraic equations for generating functions for words with $r$ occurrences of each letter and with exactly one $123$ pattern.

\vspace{0.3 in}

{\bf{1. Introduction}} ~\\

Recall that a word $w=w_1...w_k$ is an ordered list of letters on some alphabet. To say a word contains a pattern (a certain permutation of $\{1,...,m\}$) $\sigma$ is to say there exist $1 \leq i_1  <i_2<...<i_m \leq k$ such that the subword $w_{i_1}...w_{i_m}$ is {\it order isomorphic} to $\sigma$ (for example, 246 is order isomorphic to 123).  A word avoids the pattern $\sigma$ if it does not contain it.~\\ 

For a lucid history on the study of {\it forbidden patterns}, readers are welcome to refer to the introduction of Shar and Zeilberger's paper [SZ]. We also found the organization of [SZ] to be nice and suitable for the development of our paper, and we will be more or less following that. ~\\

In the present article, we say that a word $w$ in the alphabet $\{a_1,a_2,...,a_n\}$ $(a_1<a_2<...<a_n)$ is {\it{associated}} with the list $[l_1,...,l_n]$ if $w$ has $l_i$ many $a_i$'s in it, for $i$ from $1$ to $n$. For example, 231113233 is a word associated with the list $[3,2,4]$, and 223344 is a word associated with the list $[2,2,2]$. Without specifying, our default alphabet will be $\{1,...,n\}$ for some $n \geq 1$.  ~\\

In the second section, we will generalize Zeilberger's bijective proof [Z1] (a shortened version of Alexander Burstein's elegant combinatorial proof [Bu]) that the number of permutations of $\{1,...,n\}$ that contain the pattern $321$ exactly once equals $\frac{3}{n} \binom{2n}{n+3}$ and apply it to words. Although no closed form formula was found, we have a summation whose summands are expressions involving enumeration of $123$-avoiding words (for details, see Theorem $1$).~\\

In the third section, we will study, using ideas from the second section, how to extend Shar and Zeilberger's work [SZ] on generating functions enumerating $123$-avoiding words (with $r$ occurrences of each letter) to words (with $r$ occurrences of each letter) having exactly one pattern $123$. More precisely, for every positive integer $r$, Shar and Zeilberger found an algorithm for finding the defining algebraic equation for the ordinary generating function enumerating $123$-avoiding words of length $rn$ where each of the $n$ letters of $\{1,2,...,n\}$ occurs exactly $r$ times. ~\\

We will present an algorithm for finding an analogue of that, that is, a defining algebraic equation for the ordinary generating function enumerating words of length $rn$ where each of the $n$ letters of $\{1,2,...,n\}$ occurs exactly $r$ times, now with exactly one pattern $123$. We used the same (as in Shar and Zeilberger's paper [SZ]) memory-intensive, and exponential time, Buchberger's algorithm for finding Gr\"obner bases, and our computer (running Maple) found the defining algebraic equation for $r=2$: 
$${x}^{4} \left( x+4 \right) ^{2}{F}^{4}+2\,{x}^{3} \left( x+4 \right) 
 \left( 11\,x+23 \right) {F}^{3}-4\,x \left( 3\,{x}^{4}-10\,{x}^{3}-97
\,{x}^{2}-146\,x+1 \right) {F}^{2}
$$
$$
+ \left( -168\,{x}^{4}-840\,{x}^{3}-
744\,{x}^{2}+336\,x-24 \right) F+144\,{x}^{3} \left( x+2 \right)=0. 
$$

This took about a second. The minimal algebraic equation for $r=3$ has 12 as the highest power for $F$ and the computation took about $20$ seconds. Interested readers can find it on the website accompanying this paper: 
{\bf \url {http://sites.math.rutgers.edu/~my237/One123}}. The case when $r=4$ already took too long to compute (more than a month). ~\\

Now, let $a_r(n)$ be the number words of length $rn$ where each of the $n$ letters of $\{1,2,...,n\}$ occurs exactly $r$ times, with exactly one pattern $123$. In the last section of the present article, we will use the Maple package {\bf SCHUTZENBERGER}  to derive recurrence relations for our sequences. Having obtained the defining algebraic equations of the generating functions for $a_r(n)$ in the cases $r=2$ and $r=3$, Manuel Kauers kindly helped us in finding the asymptotics for our sequences $a_2(n)$ and  $a_3(n)$ (thanks to Kauers, the constants in front are fully rigorous and were computed via a step by step procedure; for details, please refer to [KP]):

$$a_2(n)=\frac{3(13-\sqrt{21})}{49}  \cdot \frac{1}{\sqrt{\pi}} \cdot 12^n \cdot n^{-3/2} \cdot (1+O(n^{-1})), $$

$$a_3(n)=\frac{-7+6\sqrt{7}}{56} \cdot \frac{1}{\sqrt{\pi}} \cdot 32^n \cdot n^{-3/2} \cdot (1+O(n^{-1})).
$$

We noticed how similar these are to the asymptotics of the sequences enumerating $123$-avoiding words with $r$ occcurences of each letter, given on page $8$ of [SZ], and we have a similar conjecture as on page $3$ of [SZ] that $a_r(n)$ is asymptotically $C_r\cdot ((r+1)2^r)^n \cdot n^{-3/2}$, where $C_r$ is a constant depending on $r$ (possibly $\frac{1}{\sqrt{\pi}}$ times a fraction of expressions involving square roots). 

\vspace{0.2 in}
{\bf 2.1. Zeilberger's shortened version of Burstein's proof on permutations containing 321 exactly once} 
\vspace{0.16 in}

In a paper published in 2011 [Bu], Burstein gave an elegant combinatorial proof of John Noonan's theorem [N] that the number of permutations of $\{1,...,n\}$ that contain the pattern $321$ exactly once equals $\frac{3}{n} \binom{2n}{n+3}$. Zeilberger [Z1] was able to shorten Burstein's proof by using a bijection between a permutation with exactly one pattern $321$, denoted as $\pi_1c\pi_2b\pi_3a\pi_4$ ($a<b<c$), with the pair $(\pi_1b\pi_2a, c\pi_3b\pi_4)$ where $\pi_1b\pi_2a$ is a $321$-avoiding permutation of $\{1,...,b\}$ and $c\pi_3b\pi_4$ is a $321$-avoiding permutation of $\{b,...,n\}$. Readers are encouraged to read Zeilberger's proof as a motivation and warm-up. Below we will see how we can use the same logic and apply it to words. ~\\

{\bf 2.2. Extension to words} ~\\

{\bf Theorem 1.}
Let $A(l_1,...,l_n)$ be the number of $123$-avoiding words associated with the list $[l_1,...,l_n]$. Then the number of words associated with list $[l_1,...,l_n]$ that contain the pattern 123 exactly once is: 
$$\sum_{b=2}^{n-1} {\sum_{j=0}^{l_{b}-1}(A(l_1,...,l_{b-1},j+1)-A(l_1,...,l_{b-1},j)) \cdot (A(l_b-j,l_{b+1},...,l_n)-A(l_b-j-1,l_{b+1},...,l_n))}.$$
\vspace{0.1in}

Before we start the proof, let us first define what a {\it good} pair of words is. Fix any $2\leq b \leq n-1$ and a list $[l_1,...,l_n]$. For any $0\leq j \leq l_b-1$ , the pair of words $(\sigma_1,\sigma_2)$ is {\it good} if $\sigma_1$ is a 123-avoiding word in $\{1,...,b\}$ associated with the list $[l_1,...,l_{b-1},j+1]$ that does not start with $b$ and $\sigma_2$ is a 123-avoiding word in $\{b,...,n\}$ associated with the list $[l_b-j,l_{b+1},...,l_n]$ that does not end with $b$. We will also say $\sigma_i$ ($i=1,2$) is good if it belongs to a good pair $(\sigma_1,\sigma_2)$.~\\

{\it Proof}. Any word $w$ associated with the list $[l_1, l_2,..., l_n]$ with exactly one pattern $123$ can be written as $\pi_1a\pi_2b\pi_3c\pi_4$ ($a<b<c$), where $abc$ is the unique $123$ pattern. All entries to the left of $b$, except $a$, must be greater than or equal to $b$, and all the entries to the right of $b$, except for $c$, must be smaller than or equal to $b$. Also, $\pi_2$ and $\pi_3$ must not contain any $b$'s, otherwise there will be another $123$ pattern. Observe that $a\pi_3b\pi_4$ is a word in $\{1,2,...,b\}$ avoiding pattern $123$ and does not start with $b$ and $\pi_1b\pi_2c$ is a word in $\{b,b+1,...,n\}$ avoiding pattern $123$ and does not end with $b$. Therefore $(a\pi_3b\pi_4,\pi_1b\pi_2c)$ is a good pair. ~\\

We now verify that there is indeed a bijection from the set of words having exactly one pattern $123$ to the set of good pairs $(\sigma_1,\sigma_2)$ (for $2\leq b\leq n-1$ and $0 \leq j \leq l_b-1$).~\\

Fix $b$ and $j$ ($2\leq b\leq n-1$, $0 \leq j \leq l_b-1$). Given a word $\pi_1a\pi_2b\pi_3c\pi_4$ (that has exactly one $123$ pattern: $abc$), we can easily map it to a unique good pair $(a\pi_3b\pi_4,\pi_1b\pi_2c)$ by first finding out what $a$, $c$ are. This is easy since we have only one $123$ pattern. Conversely, given a good pair $(\sigma_1,\sigma_2)$ , we take the first letter of $\sigma_1$ as ``$a$'' and the leftmost occurrence of $b$ as ``$b$'' and get $\pi_3$ and $\pi_4$ ($\sigma_1=a\pi_3b\pi_4$). Similarly, we take the last letter of $\sigma_2$ as ``$c$'' and the rightmost occurrence of $b$ as ``$b$'' and get $\pi_1$ and $\pi_2$ ($\sigma_2=\pi_1b\pi_2c$). Putting everything together we get a unique $\pi_1a\pi_2b\pi_3c\pi_4$. ~\\

Now, for any $b$ and $j$, the number of good $\sigma_1$ is $A(l_1,...,l_{b-1},j+1)-A(l_1,...,l_{b-1},j)$ and the number of good $\sigma_2$ is $A(l_b-j,l_{b+1},...,l_n)-A(l_b-j-1,l_{b+1},...,l_n)$. Therefore the number of words $\pi_1a\pi_2b\pi_3c\pi_4$ with exactly one pattern $123$ is: $(A(l_1,...,l_{b-1},j+1)-A(l_1,...,l_{b-1},j)) \cdot (A(l_b-j,l_{b+1},...,l_n)-A(l_b-j-1,l_{b+1},...,l_n))$. Summing over all $b$ and $j$, we get the desired result.~\\

{\bf Corollary 1.} The number of words associated with list $[l_1,...,l_n]$ that contain exactly one pattern 123 is equal to the number of words associated with list $[l_n,...,l_1]$ that contain exactly one pattern 123.  ~\\

{\it Proof.} By Theorem 1, the number of words with exactly one pattern 123 and associated with list $[l_1,...,l_n]$ is 
$$\sum_{b=2}^{n-1} {\sum_{j=0}^{l_{b}-1}(A(l_1,...,l_{b-1},j+1)-A(l_1,...,l_{b-1},j))}$$
$$\hspace{1.9 in}\cdot (A(l_b-j,l_{b+1},...,l_n)-A(l_b-j-1,l_{b+1},...,l_n)) \hspace{1.32in} (1)$$ 

\vspace{0.1in}

and the number of words with exactly one pattern 123 and associated with list $[l_n,...,l_1]$ is 
$$\sum_{b=2}^{n-1} {\sum_{j=0}^{l_{n-b+1}-1}(A(l_n,...,l_{n-b+2},j+1)-A(l_n,...,l_{n-b+2},j))}$$ 
$$\hspace{1.9in}\cdot (A(l_{n-b+1}-j,l_{n-b},...,l_1)-A(l_{n-b+1}-j-1,l_{n-b},...,l_1)). \hspace{0.65in} (2)$$

\vspace{0.1in}

When $b=k$ $(2\leq k \leq n-1)$, the inner sum of (1) becomes
$$ {\sum_{j=0}^{l_{k}-1}(A(l_1,...,l_{k-1},j+1)-A(l_1,...,l_{k-1},j))}$$
$$\hspace{1.8 in} \cdot (A(l_k-j,l_{k+1},...,l_n)-A(l_k-j-1,l_{k+1},...,l_n)) \hspace{1.38in} (3)$$ 

while when $b=n-k+1$ $(2 \leq k \leq n-1$, notice this is the ``symmetric counterpart'' of $b=k)$, the inner sum of (2) becomes
$$ {\sum_{j=0}^{l_{k}-1}(A(l_n,...,l_{k+1},j+1)-A(l_n,...,l_{k+1},j))}$$
$$\hspace{1.8 in}\cdot (A(l_k-j,l_{k-1},...,l_1)-A(l_k-j-1,l_{k-1},...,l_1)). \hspace{1.35in} (4)$$ 

\vspace{0.1in}

We only need to show $(3)=(4)$ in order to show $(1)=(2)$. Now notice that when $j=t$ $(0 \leq t \leq l_k-1)$, the summand of (3) is ~\\

$ (A(l_1,...,l_{k-1},t+1)-A(l_1,...,l_{k-1},t)) \cdot (A(l_k-t,l_{k+1},...,l_n)-A(l_k-t-1,l_{k+1},...,l_n)) \quad \quad (5)$~\\

and when $j=l_{k}-1-t$ (the ``symmetric counterpart'' of $j=t$), the summand of (4) is ~\\

$ (A(l_n,...,l_{k+1},l_k-t)-A(l_n,...,l_{k+1},l_k-t-1)) \cdot (A(t+1,l_{k-1},...,l_1)-A(t,l_{k-1},...,l_1)). \quad \quad (6)$~\\

After a small rearrangement, we can see $(5)=(6)$ because of an important result that $A(l_1,...,l_n)$ is symmetric in its arguments (meaning we can re-order the list $[l_1,...,l_n]$ in any order we want and get the same number; for details, see [SZ], page 4). Therefore as $j$ ranges from 0 to $l_k-1$, we have $(3)=(4)$. And as $b$ ranges from $0$ to $n-1$, we have $(1)=(2)$.

\vspace{0.3 in}

{\bf Corollary 2.} Fix a list $L:=[l_1,...,l_n]$. The number of words associated with L that contain exactly one pattern 123 is equal to the number of words associated with L that contain exactly one pattern 321.~\\

{\it Proof.} Let $S_1$ be the set of words associated with $[l_1,...,l_n]$ that contain exactly one pattern 123 and $S_2$ be the set of words associated with $[l_n,...,l_1]$ that contain exactly one pattern 321. Take any $w_1 \in S_1$, we can map it to a word $w_2$ associated with $[l_n,...,l_1]$ by mapping letter $i$ to letter $n-i+1$ (for all $i$ from 1 to $n$). For example, the word 121322 is mapped to 323122. Observe that $w_2$ must contain exactly one pattern 321, which occurs at the same location in $w_2$ as the location of the 123 pattern in $w_1$. Therefore $w_2 \in S_2$. Clearly this is a bijection from $S_1$ to $S_2$. So $|S_1|=|S_2|$. This along with Corollary 1 gives Corollary 2. ~\\

{\bf Remark.} 
\vspace{0.2in}

One may wonder if the number of words associated with $[l_1,...,l_n]$ that contain exactly one pattern 123 is equal to the number of words associated with $[l_1,...,l_n]$ that contain exactly one pattern 132. This is not the case (if this were the case, we would have an analogue of the result that the 123-avoiding words associated with $[l_1,...,l_n]$ are equinumerous with the 132-avoiding words associated with $[l_1,...,l_n]$, see [Z2]). For example, the number of permutations of $\{1,2,...,n\}$ ($n\geq 1$) that contain exactly one pattern 123 is $\frac{3}{n} \binom{2n}{n+3}$ [Z1] while the number of permutations of $\{1,2,...,n\}$ ($n\geq 1$) that contain exactly one pattern 132 is $\binom{2n-3}{n-3}$ [B\'{o}].

\newpage

{\bf 3.1. Some crucial background and generating functions for words avoiding pattern 123} ~\\

In the beautiful paper by Shar and Zeilberger [SZ], methods for finding the algebraic equation for the ordinary generating function enumerating $123$-avoiding words of length $rn$, where each of the $n$ letters of $\{1,2,...,n\}$ occurs exactly $r$ times were given. First we present some important definitions and results of that paper here.~\\

For $0 \leq i \leq j \leq r-1$ and $n \geq 0$, let $W_r^{(i,j)}(n)$ be the set of $123$-avoiding words of length $rn+i+j$, in the alphabet $\{1,2,..n,n+1,n+2\}$, with $i$ occurrences of the letter $1$, $j$ occurrences of the $n+2$, and exactly $r$ occurrences of the other $n$ letters, and let $W_r^{(i,j)}$ be the union of $W_r^{(i,j)}(n)$ over all $n \geq 0$. Let $g_r^{(i,j)}(x)$ be the {\it weight enumerator} for $W_r^{(i,j)}$, with respect to the weight $w \rightarrow x^{length(w)}$.
(Note that the $W_r^{(i,j)}$'s have the same weight enumerator if any two letters have $i$ and $j$ occurrences respectively, and the remaining letters each occurs exactly $r$ times. For a detailed explanation, see [SZ], page 3-4.)~\\

Shar and Zeilberger were able to find a system of $\binom {r+1}{2}$ equations for $g_r^{(i,j)}(x)$ ($0\leq i\leq j \leq r-1$), with the convention that if $s>k$ then $g_r^{(s,k)}=g_r^{(k,s)}$: 
$$g_r^{(i,j)}(x)=\delta_{i,0}\delta_{j,0}+x\sum_{t=0}^{r-1}g_r^{(i,t)}(x)g_r^{((r-t) \ \text{mod}\ r,\ (j-1) \ \text{mod} \ r)}(x)+\sum_{m=0}^{i-1}x^{m+1}g_r^{(i-m,j-1)}(x)$$

For example, in the case when $r=2$, we would get the following system of equations: ~\\

$g_2^{(0,0)}(x)=1+xg_2^{(0,0)}(x)g_2^{(0,1)}(x)+xg_2^{(0,1)}(x)g_2^{(1,1)}(x)$ ~\\

$g_2^{(0,1)}(x)=xg_2^{(0,0)}(x)^2+xg_2^{(0,1)}(x)^2$ ~\\

$g_2^{(1,1)}(x)=xg_2^{(0,0)}(x)g_2^{(0,1)}(x)+xg_2^{(0,1)}(x)(1+g_2^{(1,1)}(x))$ ~\\

Solving this system of equations in the three unknowns $g_2^{(0,0)}(x),g_2^{(0,1)}(x),g_2^{(1,1)}(x)$, we get the weight enumerators for $W_2^{(0,0)}$, $W_2^{(0,1)}$ and $W_2^{(1,1)}$. ~\\

Once we have the weight enumerators, we can easily get the corresponding generating functions by doing a little operation. For example, if we have an explicit expression for $g_2^{(0,0)}(x)$ ($g_2^{(0,0)}(x)=1+x^2+6x^4+43x^6+352x^8+3114x^{10}+...$), the corresponding generating function will be $f_2^{(0,0)}(x)=1+x+6x^2+43x^3+352x^4+3114x^5+...$ (that is, $f_2^{(0,0)}(x)=g_2^{(0,0)}(x^{1/2}$)). ~\\

\newpage

{\bf 3.2. Extension to generating functions for words with exactly $r$ occurrences of each letter, and with exactly one pattern $123$}~\\

{\bf Definitions:} Let $V_r(n)$ be the set of words in the alphabet $\{1,...,n\}$ with exactly $r$ occurrences of each letter, and with exactly one pattern $123$. Let $V_r$ be $\bigcup_{n=0}^{\infty}V_r(n)$.~\\

Let $h_r(x)$ be the weight enumerator for $V_r$ (as always, with weight $w \rightarrow x^{length(w)})$ and let $f_r(x)$ be the corresponding generating function. We will be following the framework of [SZ], with two warm-up cases leading to the general case. ~\\

{\bf First warm-up: $r=1$} ~\\

Claim: $h_1(x)=(g_1^{(0,0)}(x)-xg_1^{(0,0)}(x)-1)^2/x$. ~\\

{\it Proof}. Recall that $g_1^{(0,0)}(x)$ $(=1+x+2x^2+5x^3+14x^4+42x^5+...)$ is the weight enumerator for $123$-avoiding permutations on $\{1,...,n,...\}$. We prove this claim by showing that the coefficient of $x^{n}$ ($n\geq0$) on the right hand side is exactly the number of good pairs $(a\pi_3b\pi_4,\pi_1b\pi_2c)$ ($2 \leq b \leq n-1 $), which equals to the number of permutations on $\{1,...,n\}$ with exactly one pattern $123$ (by Zeilberger's proof [Z1]).  ~\\

For any fixed $b$ ($2 \leq b \leq n-1$), a good $a\pi_3b\pi_4$ would be a $123$-avoiding permutation on $\{1,...,b\}$ that does not start with $b$. Similarly, a good $\pi_1b\pi_2c$ would be a $123$-avoiding permutation on $\{b,...,n\}$ that does not end with $b$.~\\

Note that the coefficient of $x^{b}$ in $g_1^{(0,0)}(x)-xg_1^{(0,0)}(x)-1$ is exactly the number of good $a\pi_3b\pi_4$. 
(The $x$ in front of $g_1^{(0,0)}(x)$ corresponds to having $b$ in front of a permutation, and the $-1$ corresponds to an empty permutation. We don't want either of these.) ~\\

Similarly, the coefficient of $x^{n-b+1}$ in $g_1^{(0,0)}(x)-xg_1^{(0,0)}(x)-1$ is the number of good $\pi_1b\pi_2c$. Multiplying the two, we have that the coefficient of $x^{n+1}$ ($=x^{b} \cdot x^{n-b+1}$) in 

\begin{center}
$(g_1^{(0,0)}(x)-xg_1^{(0,0)}(x)-1)^2$
\end{center}
is the number of good pairs $(a\pi_3b\pi_4,\pi_1b\pi_2c)$ ($b$ ranges from $2$ to $n-1$). Dividing by $x$, we get the coefficient of $x^{n}$ in
\begin{center}
$(g_1^{(0,0)}(x)-xg_1^{(0,0)}(x)-1)^2/x$ ~\\
\end{center}
is the number of good pairs $(a\pi_3b\pi_4,\pi_1b\pi_2c).$~\\

\vspace{0.1in}

{\bf Second warm-up: $r=2$} ~\\

Claim: $h_2(x)=2 \cdot (g_2^{(0,0)}(x)-xg_2^{(0,1)}(x)-1)(g_2^{(0,1)}(x)-xg_2^{(0,0)}(x))/x$. ~\\

{\it Proof}. Recall that $g_2^{(0,0)}(x)$ ($=1+x^2+6x^4+43x^6+...$) is the weight enumerator for $123$-avoiding words on $\{1,1,...,n,n,...\}$ (or equivalently, $123$-avoiding words associated with $[2,2,...]$) and $g_2^{(0,1)}(x)$ $(=x+3x^3+19x^5+145x^7...)$ is the weight enumerator for $123$-avoiding words on $\{1,2,2,3,3,...,n,n,...\}$. As in the first warm-up, we prove this claim by showing that the coefficient of $x^{2n}$ ($n\geq 0$) on the right hand side is exactly the number of good pairs $(a\pi_3b\pi_4,\pi_1b\pi_2c)$ ($2 \leq b \leq n-1 $), which equals to the number of words on $\{1,1,...,n,n\}$ with exactly one pattern $123$ (by the proof of Theorem $1$).  ~\\

For any $b$ ($2 \leq b \leq n-1$), we have the following two cases: either $\pi_4$ contains one $b$ and $\pi_1$ contains no $b$ or the other way around.~\\

{\bf Case 1}: $\pi_4$ contains one $b$ and $\pi_1$ contains no $b$.~\\

Then a good $a\pi_3b\pi_4$ would be a $123$-avoiding word on $\{1,1,...,b,b\}$ that does not start with $b$. Similarly, a good $\pi_1b\pi_2c$ would be a $123$-avoiding word on $\{b,b+1,b+1,...,n,n\}$ that does not end with $b$.~\\

Note that the coefficient of $x^{2b}$ in $g_2^{(0,0)}(x)-xg_2^{(0,1)}(x)-1$ is exactly the number of good $a\pi_3b\pi_4$. (The $x$ in front of $g_2^{(0,1)}(x)$ corresponds to having $b$ in front of a word, and the $-1$ corresponds to an empty word. We don't want either of these.) ~\\

Similarly, the coefficient of $x^{2(n-b)+1}$ in $g_2^{(0,1)}(x)-xg_2^{(0,0)}(x)$ is the number of good $\pi_1b\pi_2c$.  ~\\

Multiplying the two, and let $b$ range from 2 to $n-1$, we have that the coefficient of $x^{2n+1}$ ($=x^{2b} \cdot x^{2(n-b)+1}$) in 
\begin{center}
$(g_2^{(0,0)}(x)-xg_2^{(0,1)}(x)-1)(g_2^{(0,1)}(x)-xg_2^{(0,0)}(x))$ 
\end{center}
is the number of good pairs $(a\pi_3b\pi_4,\pi_1b\pi_2c)$ if the additional $b$ is in $\pi_4$. Dividing by $x$, we get the coefficient of $x^{2n}$ in
\begin{center}
$(g_2^{(0,0)}(x)-xg_2^{(0,1)}(x)-1)(g_2^{(0,1)}(x)-xg_2^{(0,0)}(x))/x$ 
\end{center}
is the number of good pairs $(a\pi_3b\pi_4,\pi_1b\pi_2c)$ if the additional $b$ is in $\pi_4$. ~\\

{\bf Case 2}: $\pi_1$ contains one $b$ and $\pi_4$ contains no $b$. ~\\

Then a good $a\pi_3b\pi_4$ would be a $123$-avoiding word on $\{1,1,...,b-1,b-1,b\}$ that does not start with $b$. A good $\pi_1b\pi_2c$ would be a $123$-avoiding word on $\{b,b,...,n,n\}$ that does not end with $b$.~\\

Now, the coefficient of $x^{2b-1}$ in $g_2^{(0,1)}(x)-xg_2^{(0,0)}(x)$ is exactly the number of good $a\pi_3b\pi_4$. Similarly, the coefficient of $x^{2(n-b)+2}$ in $g_2^{(0,0)}(x)-xg_2^{(0,1)}(x)-1$ is the number of good $\pi_1b\pi_2c$.  ~\\

Multiplying the two, we have that the coefficient of $x^{2n+1}$ ($=x^{2b-1} \cdot x^{2(n-b)+2}$) in 
\begin{center}
$(g_2^{(0,1)}(x)-xg_2^{(0,0)}(x))(g_2^{(0,0)}(x)-xg_2^{(0,1)}(x)-1)$ 
\end{center}
is the number of good pairs $(a\pi_3b\pi_4,\pi_1b\pi_2c)$ ($b$ ranges from 2 to $n-1$) if the additional $b$ is in $\pi_1$. Dividing by $x$, we get the coefficient of $x^{2n}$ in 
\begin{center}
$(g_2^{(0,1)}(x)-xg_2^{(0,0)}(x))(g_2^{(0,0)}(x)-xg_2^{(0,1)}(x)-1)/x$ 
\end{center}
is the number of good pairs $(a\pi_3b\pi_4,\pi_1b\pi_2c)$ if the additional $b$ is in $\pi_1$. ~\\

Therefore the coefficient of $x^{2n}$ in $2 \cdot (g_2^{(0,0)}(x)-xg_2^{(0,1)}(x)-1)(g_2^{(0,1)}(x)-xg_2^{(0,0)}(x))/x$ is the number of good pairs $(a\pi_3b\pi_4,\pi_1b\pi_2c)$, which is equal to the number of words in $\{1,1,...,n,n\}$ that have exactly one pattern $123$. That is, we have shown the weight enumerator for $V_2$ is as claimed to be. To get the generating function $f_2(x)$ for $V_2$ we simply let $f_2(x)=h_2(x^{1/2})$. ~\\

Readers are welcome to compare $h_2(x)$ with the earlier formula in the case when $l_i=2 (1 \leq i \leq n)$: 
$$\sum_{b=2}^{n-1} {\sum_{j=0}^{1}(\underbrace{A(j+1,2,2,...,2)}_\text{$b-1$ many 2's}-\underbrace{A(j,2,2,...,2)}_\text{$b-1$ many 2's}) \cdot (\underbrace{A(2-j,2,2,...,2)}_\text{$n-b$ many 2's}-\underbrace{A(1-j,2,2,...,2)}_\text{$n-b$ many 2's})}.$$~\\

{\bf The general case} ~\\

{\bf Theorem 2.} ~\\$$h_r(x)=\frac{1}{x}\sum_{i=1}^{r}(g_r^{(0,\ i \ \text{mod} \ r)}-xg_r^{(0,\ i-1)}-\delta_{(i \ \text{mod} \ r, \ 0)})(g_r^{(0,\ (r+1-i) \ \text{mod} \ r)}-xg_r^{(0,\ r-i)}-\delta_{((r+1-i)\ \text{mod} \ r,\ 0)}).$$~\\

The general case is derived using the exact same logic as for the warm-up cases. Instead of having two cases as in the second warm-up, here we have $r$ cases. Interested readers are welcome to verify the formula for $r=3$ by themselves, and the general case should be apparent after this verification. As before, to get the generating function for $V_r$ we simply let $f_r(x)=h_r(x^{1/r})$.

\newpage

{\bf 4. Using Maple packages}~\\

As noted in Shar and Zeilberger's paper ([SZ], page $7$), now that we know $f_r(x)$ has the property of being algebraic, the sequence ${a_r(n)}$ satisfies some homogeneous linear recurrence equation with polynomial coefficients.~\\

Using the {\it algtorec} procedure in the {\bf SCHUTZENBERGER} package written by Doron Zeilberger, available from: ~\\

{\bf \url {http://www.math.rutgers.edu/~zeilberg/tokhniot/SCHUTZENBERGER.txt}}~\\

we are able to find (rigorously) recurrences (in operator notation) for our sequences when $r=1$ and $r=2$ ($r=3$ took too long to compute): ~\\

For $r=1$ we get: 
$$(2\,n \left( 2\,n+1 \right) - \left( n+4 \right)  \left( n-2 \right)N)a_1(n)=0$$ ~\\
(which agrees with the already known formula $a_1(n)=\frac{3}{n} \binom{2n}{n+3}$).~\\

In the case when $r=2$, {\it algtorec} returned a operator of degree $8$, but it can be reduced to a minimal operator of degree $4$ (thanks to Manuel Kauers for pointing it out), that is: ~\\

$$(36(1 + n)(2 + n)(1 + 2n)(3 + 2n)(18154800 + 23101940n +
10635771n^2 + 2093616n^3 + 147833n^4)
$$
$$+12(2 + n)(3 + 2n)(1283329440 + 3700267618n + 4200957553n^2 +
2408049238n^3 + 735936616n^4 
$$
$$
+113774584n^5 + 6948151n^6)N + (282564806400 + 1066356868608n +
1704365727480n^2 
$$
$$+1511140337906n^3+814587362081n^4 + 273775889012n^5 + 56080140110n^6 +
6405068474n^7
$$
$$
+ 312371129n^8)N^2 -2(4 + n)(11939685120 + 40890299130n + 56943840213n^2 +
41794221496n^3 
$$
$$
+17488032270n^4 + 4183030930n^5 +531527997n^6 + 27792604n^7)N^3 + 8(1 + n)(4 + n)(5 + n)(11
$$
$$+ 2n)(3742848 + 7519914n + 5241921n^2 + 1502284n^3 + 147833n^4)N^4)a_2(n)=0$$
\vspace{0.1in}
  
Everything in this paper is implemented (with explanation) in the Maple packages {\bf Words123New} and {\bf PW123} and available from: {\bf \url {http://sites.math.rutgers.edu/~my237/One123}}, which also includes some sample input and output files. ~\\

\newpage
{\bf Acknowledgements}~\\

The author thanks Doron Zeilberger for bringing this project into her attention and continued conversation. The author also would like to thank Matthew Russell for offering suggestions that improved the exposition of this paper, Andrew Lohr and Alejandro Ginory for proofreading the draft, and Manuel Kauers for his help in finding the asymptotics.~\\

\vspace{0.2in}

{\bf References}~\\

[B\'{o}] M. B\'{o}na, {\it Permutations with one or two 132-subsequences},
Discrete Math. 181 (1998), 267-274.~\\
{\bf \url{http://www.sciencedirect.com/science/article/pii/S0012365X97000629?via%3Dihub}}~\\

[Bu] A. Burstein, {\it A short proof for the number of permutations containing pattern 321 exactly once}, The Electronic J. of Combinatorics 18(2) (2011), P21. ~\\
{\bf \url {http://www.combinatorics.org/ojs/index.php/eljc/article/view/v18i2p21/pdf}}~\\

[K] M. Kauers, A Mathematica Package for Computing Asymptotic
Expansions of Solutions of P-Finite Recurrence Equations, Technical
Report RISC 11-04, Johannes Kepler University, April 2011.~\\

[KP] M. Kauers and P. Paule, {\it The Concrete Tetrahedron: Symbolic Sums, Recurrence Equations, Generating Functions, Asymptotic Estimates},  Springer (2011), Section 6.5.~\\

[N] J. Noonan, {\it The number of permutations containing exactly one increasing subsequence of length three}, Discrete Math. 152 (1996), 307-313.~\\
{\bf \url{http://www.sciencedirect.com/science/article/pii/0012365X9500247T}}~\\

[SZ] N. Shar and D. Zeilberger, {\it The (Ordinary) Generating Functions Enumerating $123$-Avoiding Words with $r$ occurences of each of $1,2,...,n$ are Always Algebraic}, Annals of Combinatorics 20 (2016), 387-396. ~\\
{\bf \url{http://sites.math.rutgers.edu/~zeilberg/mamarim/mamarimhtml/words123.html}}~\\

[Z1] D. Zeilberger, {\it Alexander Burstein's Lovely Combinatorial Proof of John Noonan's Beautiful Theorem that the number of $n$-permutations that contain the Pattern $321$ Exactly Once Equals $(3/n)(2n!)/((n-3)!(n+3)!)$}, Personal Journal of Shalosh B. Ekhad and Doron Zeilberger, Oct.18, 2011.~\\
{\bf \url{http://sites.math.rutgers.edu/~zeilberg/mamarim/mamarimhtml/burstein.html}}~\\

[Z2] D. Zeilberger, {\it A Snappy Proof That 123-Avoiding Words are Equinumerous With 132-Avoiding Words}, Personal Journal of Shalosh B. Ekhad and Doron Zeilberger, April 11, 2005.
{\bf \url{http://sites.math.rutgers.edu/~zeilberg/mamarim/mamarimhtml/a123.html}}
\vspace{0.2in}

\par\noindent\rule{\textwidth}{0.4pt}

Department of Mathematics, Rutgers University, Piscataway, NJ 08854~\\
\it{E-mail address}: \bf{my237@math.rutgers.edu}
\end{document}